\documentclass[11pt]{article}
\usepackage{amsfonts}
\usepackage{amssymb}
\usepackage{amsmath}

\def\sqr#1#2{{\vcenter{\vbox{\hrule height.#2pt
              \hbox{\vrule width.#2pt height#1pt \kern#1pt \vrule
width.#2pt}
              \hrule height.#2pt}}}}
\def\signed #1{{\unskip\nobreak\hfil\penalty50
              \hskip2em\hbox{}\nobreak\hfil#1
              \parfillskip=0pt \finalhyphendemerits=0 \par}}
\def\endpf{\signed {$\sqr69$}}
\def\dbR{{\mathbb{R}}}

\def\g{\gamma}
\def\d{\delta}
\def\e{\varepsilon}
\def\z{\zeta}

\def\l{\lambda}
\def\m{\mu}

\def\si{\sigma}
\def\t{\tau}
\def\f{\varphi}

\def\3n{\negthinspace \negthinspace \negthinspace }
\def\2n{\negthinspace \negthinspace }
\def\1n{\negthinspace }
\def\ns{\noalign{\smallskip} }
\def\ns{\noalign{\medskip} }
\def\ds{\displaystyle}
\def\D{\Delta}
\def\Th{\Theta}

\def\F{\Phi}

\def\cF{{\cal F}}

\def\cS{{\cal S}}

\def\cU{{\cal U}}

\def\no{\noindent}

\def\ms{\medskip}
\def\bs{\bigskip}
\def\q{\quad}
\def\qq{\qquad}
\def\hb{\hbox}

\def\Ra{\mathop{\Rightarrow}}

\def\lan{\mathop{\langle}}
\def\ran{\mathop{\rangle}}

\def\pa{\partial}

\def\wt{\widetilde}
\def\cd{\cdot}
\def\cds{\cdots}

\def\({\Big (}
\def\){\Big )}
\def\[{\Big[}
\def\]{\Big]}

\def\be{\begin{equation}}
\def\bel{\begin{equation}\label}
\def\ee{\end{equation}}
\def\bea{\begin{eqnarray}}
\def\eea{\end{eqnarray}}
\def\bt{\begin{theorem}}
\def\et{\end{theorem}}
\def\bc{\begin{corollary}}
\def\ec{\end{corollary}}
\def\bl{\begin{lemma}}
\def\el{\end{lemma}}
\def\bp{\begin{proposition}}
\def\ep{\end{proposition}}
\def\br{\begin{remark}}
\def\er{\end{remark}}
\def\ba{\begin{array}}
\def\ea{\end{array}}
\def\bd{\begin{definition}}
\def\ed{\end{definition}}

\newtheorem{lemma}{Lemma}[section]
\newtheorem{remark}{Remark}[section]

\newtheorem{theorem}{Theorem}[section]
\newtheorem{corollary}{Corollary}[section]

\newtheorem{definition}{Definition}[section]
\newtheorem{proposition}{Proposition}[section]

\makeatletter
   
   \@addtoreset{equation}{section}
\makeatother

\begin{document}
\title{\bf A Deterministic Linear Quadratic\\ Time-Inconsistent
Optimal Control Problem\footnote{This work is supported in part by
the NSF grant DMS-1007514.}}

\author{Jiongmin Yong\\
Department of Mathematics, University of Central Florida, Orlando,
FL 32816}

\maketitle

\begin{abstract}

A time-inconsistent optimal control problem is formulated and
studied for a controlled linear ordinary differential equation with
quadratic cost functional. A notion of equilibrium control is
introduced, which can be regarded as a time-consistent solution to
the original time-inconsistent problem. Under certain conditions, we
constructively prove the existence of such an equilibrium control
which is represented via a forward ordinary differential equation
coupled with a backward Riccati--Volterra integral equation. Our
constructive approach is based on the introduction of a family of
$N$-person non-cooperative differential games.

\end{abstract}

\ms


\ms


\section{Introduction --- Time-Consistency Issue.}\label{1}
We begin with a classical optimal control problem for an ordinary
differential equation (ODE, for short). Let $T>0$. For any {\it
initial pair} $(t,x)\in[0,T)\times\dbR^n$, consider the following
controlled ODE:
\bel{ODE1}\left\{\ba{ll}
\ns\ds\dot X(s)=b(s,X(s),u(s)),\qq s\in[t,T],\\
\ns\ds X(t)=x,\ea\right.\ee
where $b:[0,T]\times\dbR^n\times U\to\dbR^n$ is a given map,
$u(\cd)$, a function valued in some metric space $U$, is called a
{\it control}, and $X(\cd)$ is called the {\it state trajectory}. We
denote
\bel{}\cU[t,s]=\Big\{u:[t,s]\to U\bigm|u(\cd)\hb{ is
measurable}\Big\},\qq\forall0\le t\le s\le T.\ee
Under some mild conditions, for any {\it initial pair}
$(t,x)\in[0,T)\times\dbR^n$, and $u(\cd)\in\cU[t,T]$, (\ref{ODE1})
admits a unique solution $X(\cd)\equiv X(\cd\,;t,x,u(\cd))$. Then we
can introduce the following cost functional which measures the
performance of the control $u(\cd)$:
\bel{cost1}J(t,x;u(\cd))=\int_t^Tg(s,X(s;t,x,u(\cd)),u(s))ds+h(X(T;t,x,u(\cd))),\ee
for some given maps $g:[0,T]\times\dbR^n\times U\to\dbR$ and
$h:\dbR^n\to\dbR$. The terms on the right hand side of (\ref{cost1})
are referred to as the {\it running cost} and the {\it terminal
cost}, respectively. The following is a classical optimal control
problem.

\ms

\bf Problem (D). \rm For any given initial pair
$(t,x)\in[0,T)\times\dbR^n$, find a $\bar u(\cd)\in\cU[t,T]$ such
that
\bel{}J(t,x;\bar u(\cd))=\inf_{u(\cd)\in\cU[t,T]}J(t,x;u(\cd)).\ee
Any $\bar u(\cd)\in\cU[t,T]$ satisfying the above is called an {\it
optimal control} for $(t,x)$, $\bar X(\cd)\equiv X(\cd\,;t,x,\bar
u(\cd))$ is called the corresponding {\it optimal trajectory}, and
$(\bar X(\cd),\bar u(\cd))$ is referred to as an {\it optimal pair}.

\ms

Note that sometimes we might encounter the following seemingly a
little more general form of the cost functional
\bel{cost1.5}\wt
J(t,x;u(\cd))=\int_t^Te^{-\int_t^sc(r,X(r),u(r))dr}g(s,X(s),u(s))ds
+e^{-\int_t^Tc(r,X(r),u(r))dr}h(X(T)),\ee
with $c(\cd)$ being some map taking nonnegative values, which may be
called a {\it discount map}. A special case is $c(\cd)=\d>0$, a
positive constant (which is call a {\it discount rate}). Due to its
form, the term $\ds e^{-\int_t^sc(r,X(r),u(r))dr}$ is called an {\it
exponential discounting}. If we introduce
\bel{}\left\{\ba{ll}
\ns\ds\dot X^0(s)=c(s,X(s),u(s)),\qq s\in[t,T],\\
\ns\ds X^0(t)=0,\ea\right.\ee
and regard $X^0(\cd)$ an additional component of the state, then the
state equation is augmented by one dimension and the cost functional
becomes
\bel{cost1.6}\wt J(t,x;u(\cd))=\int_t^Te^{-X^0(s)}g(s,X(s),u(s))ds
+e^{-X^0(T)}h(X(T)).\ee
which is of form (\ref{cost1}). Therefore, an optimal control
problem with an exponential discounting can be transformed to an
optimal control problem without exponential discounting. In another
word, containing an exponential discounting in the cost functional
does not make the original problem mathematically more general.

\ms

Dynamic programming method is a powerful classical approach to
Problem (D). This method suggests us define the {\it value function}
of Problem (D) by the following:
\bel{value}\left\{\ba{ll}
\ns\ds
V(t,x)=\inf_{u(\cd)\in\,\cU[t,T]}J(t,x;u(\cd)),\qq(t,x)\in[0,T)\times\dbR^n,\\
\ns\ds V(T,x)=h(x),\qq\forall x\in\dbR^n.\ea\right.\ee
It is well-known that the following Bellman's {\it principle of
optimality} holds (\cite{YZ1999}):
\bel{Optimality}\ba{ll}
\ns\ds V(t,x)=\inf_{u(\cd)\in\,\cU[t,\t]}\Big\{\int_t^\t
g(s,X(s),u(s))ds+V(\t,X(\t))\Big\},\qq\forall\t\in[t,T].\ea\ee
Now, if $\bar u(\cd)\in\cU[t,T]$ is an optimal control for the
initial pair $(t,x)$, then, from the above, for any $\t\in(t,T)$,
\bel{}\ba{ll}
\ns\ds V(t,x)=J(t,x;\bar u(\cd))=\int_t^Tg(s,\bar X(s),\bar u(s))ds+h(\bar X(T))\\
\ns\ds\qq\q=\int_t^\t g(s,\bar X(s),\bar u(s))ds+J(\t,\bar
X(\t);\bar
u\big|_{[\t,T]}(\cd))\\
\ns\ds\qq\q\ge\int_t^\t g(s,\bar X(s),\bar u(s))ds+V(\t,\bar
X(\t))\ge V(t,x).\ea\ee
Hence, all the equalities in the above have to hold. Consequently,
\bel{}\inf_{u(\cd)\in\cU[\t,T]}J(\t,\bar X(\t);u(\cd))=V(\t,\bar
X(\t))=J(\t,\bar X(\t);\bar u\big|_{[\t,T]}(\cd)).\ee
This means that for any $0\le t<\t<T$, the restriction $\bar
u\big|_{[\t,T]}(\cd)\in\cU[\t,T]$ of optimal control $\bar
u(\cd)\in\cU[t,T]$ for the initial pair $(t,x)$ on the time interval
$[\t,T]$ is optimal for the initial pair $(\t,\bar X(\t))$. Such a
phenomenon is referred to as the {\it time-consistency} of Problem
(D). The advantage of the time-consistency is that one needs only to
solve Problem (D) for a given initial pair
$(t,x)\in[0,T)\times\dbR^n$, and as time goes by, the restriction of
the optimal control $\bar u(\cd)$ for $(t,x)$ on any later time
interval $[\t,T]$ will automatically be an optimal control for the
corresponding initial pair $(\t,\bar X(\t))$.

\ms

However, common sense tells us that the time-consistency issue in
real life is actually never so simple. There are two main reasons:
First, as time goes by, the environment (in the broad sense) is
changing, for example, invention of new technology, new limits of
resource allocation, etc., and therefore the controlled system has
to be modified according to the new initial pairs; and secondly,
people keep changing minds/objectives, which leads to the change of
cost functional. Due to these changes, one expects some dramatic
changes in the formulation of optimal control problems, as well as
the solutions to the problems.

\ms

To make our statement more appealing from mathematical point of
view, let us look at a very simple illustrative example. Consider a
one-dimensional controlled ODE:
\bel{1.9}\left\{\ba{ll}
\ns\ds\dot X(s)=u(s),\qq s\in[t,T],\\
\ns\ds X(t)=x,\ea\right.\ee
with cost functional
\bel{1.10}J(t,x;u(\cd))=\int_t^Tu(s)^2ds+h(t)X(T;t,x,u(\cd))^2,\ee
where $h:[0,T]\to[\d,\infty)$, for some $\d>0$, and $U=\dbR$. We
pose the following optimal control problem.

\ms

\bf Problem (C). \rm For given $(t,x)\in[0,T)\times\dbR$, find a
$\bar u(\cd)\in\cU[t,T]$ such that
\bel{1.3}J(t,x;\bar
u(\cd))=\inf_{u(\cd)\in\cU[t,T]}J(t,x;u(\cd)).\ee

\ms

Note that the above problem looks like a simple standard linear
quadratic optimal control problem (LQ problem, for short), except
that the terminal weight $h(t)$ depends on the parameter $t$ (which
is the initial time of the problem).

\ms

It is clear that for any initial pair $(t,x)\in[0,T)\times\dbR$,
$u(\cd)\mapsto J(t,x;u(\cd))$ is convex and coercive. Thus, there
exists a unique optimal control for Problem (C). We can show that
(see the Appendix) the optimal control of Problem (C) is given by
\bel{}\bar u(s)\equiv\bar u(s;t,x)=-{xh(t)\over1+h(t)(T-t)}\,,\qq
s\in[t,T],\ee
and the corresponding optimal trajectory is given by
\bel{}\bar X(s;t,x,\bar u(\cd))=x{1+h(t)(T-s)\over1+h(t)(T-t)}\,,\qq
s\in[t,T].\ee
Now, for $\t\in(t,T)$, we consider Problem (C) on $[\t,T]$ with
initial state
\bel{}y=\bar X(\t;t,x,\bar
u(\cd))=x{1+h(t)(T-\t)\over1+h(t)(T-t)}\,.\ee
We can show that
\bel{}J(\t,y;\bar u(\cd))-\inf_{u(\cd)\in\cU[\t,T]}J(\t,y;u(\cd))
={x^2[h(\t)-h(t)]^2(T-\t)\over[1+h(\t)(T-t)][1+h(t)(T-\t)]^2}\,.\ee
This tells us that the restriction of $\bar u(\cd\,;t,x)$ on
$[\t,T]$ is not optimal for Problem (C) with initial pair $(\t,\bar
X(\t;t,x))$, in general. Such a phenomenon is called {\it
time-inconsistency}.

\ms

Qualitative analysis on time-inconsistent behaviors can at least be
traced back to the works by Hume \cite{Hume} in 1739 and by Smith
\cite{Smith} in 1759. Later relevant works were made by Malthus
\cite{Malthus} in 1828, Jevons \cite{Jevons} in 1871, Marshall
\cite{Marshall} in 1890, B\"ohm--Bawerk \cite{BB} in 1891, and
Pareto \cite{Pareto} in 1909, and so on. Mathematical formulation of
time-inconsistency was firstly presented by Strotz \cite{Strotz} in
1955, followed by Pollak \cite{Pollak 1968}, Peleg--Yaari
\cite{Peleg-Yaari 1973}, Goldman \cite{Goldman 1980}, Laibson
\cite{L1997}, etc. See Palacios-Huerta \cite{Palacios-Huerta} for an
interesting survey on the history. The above-mentioned mathematical
works, starting from Strotz, mainly studied problems for either
discrete dynamic systems or simple ODEs, involving non-exponential
discounting, meaning that in the cost functional (see
(\ref{cost1.5}) with $c(\cd)=\d$), the classical exponential
discounting $e^{-\d(s-t)}$ is replaced by a function $h(s-t)$.
Recently, Ekeland--Lazrak \cite{Ekeland-Lazrak 2006} and
Ekeland--Pirvu \cite{Ekeland-Privu 2007} continued the study of
non-exponential discounting problems both for simple ODEs and SDEs.
At the same time, Basak--Chabakauri \cite{Basak-Chabakarui 2008} and
Bj\"ork--Murgoci \cite{BM2008} started to discuss the problems with
the cost/payoff functional depending on the initial pair $(t,x)$. We
refer to \cite{Grenadier-Wang2005}, \cite{Herings-Rohde},
\cite{Krusell-Smith 2003}, \cite{Miller-Salmon 1985},
\cite{Tesfatsion 1986} for some relevant results.

\ms

In general, for any given initial pair $(t,x)\in[0,T)\times\dbR^n$,
we can consider the following controlled system:
\bel{1.18}\left\{\ba{ll}
\ns\ds\dot X(s)=b(t,x,s,X(s),u(s)),\qq s\in[t,T],\\
\ns\ds X(t)=x,\ea\right.\ee
with the cost functional:
\bel{1.19}\ba{ll}
\ns\ds J(t,x;u(\cd))=\int_t^Tg(t,x,s,X(s),u(s))ds+h(t,x,X(T)).\ea\ee
We point out that state equation (\ref{1.18}) and cost functional
(\ref{1.19}) are significantly different from (\ref{ODE1}) and
(\ref{cost1}), respectively, due to the way they depend on the
initial pair $(t,x)\in[0,T)\times\dbR^n$. Such a dependence allows
us to catch some situations that people will modify the control
system and/or the cost functional at different initial pair.
Clearly, our setting is much more general than \cite{Ekeland-Lazrak
2006}. Naturally, one could pose the following optimal control
problem.

\ms

\bf Problem (N). \rm For $(t,x)\in[0,T)\times\dbR^n$, find $\bar
u(\cd)\in\cU[t,T]$ such that
\bel{}J(t,x;\bar
u(\cd))=\inf_{u(\cd)\in\,\cU[t,T]}J(t,x;u(\cd)).\ee

\ms

It is clear that Problem (C) is a special case of Problem (N).
Hence, Problem (N) is time-inconsistent, in general. Any optimal
control $\bar u(\cd)\in\cU[t,T]$ of Problem (N) is referred to as a
{\it pre-committed} optimal control on $[t,T]$. Due to the
time-inconsistency, finding an optimal control $\bar
u(\cd)\in\cU[t,T]$ for Problem (N) (assuming it exists) might not be
very useful (if it is not useless) in long run. Hence, Problem (N)
is natural, but is a little too naive.

\ms

In this paper, we will concentrate on a linear-quadratic
time-inconsistent control problem. We will present a time-consistent
solution via a ``sophisticated'' approach. The main idea comes from
the works \cite{Strotz}, \cite{Pollak 1968}, \cite{Peleg-Yaari
1973}, and \cite{Goldman 1980}. Here is a brief description. Take a
partition $\D:0=t_0<t_1<\cds<t_N=T$ of the time interval $[0,T]$.
Consider an $N$-person non-cooperative differential game: the $k$-th
player (which may be called self--$k$) starts the game from the
initial pair $(t_{k-1},X(t_{k-1}))$ and controls the system on
$[t_{k-1},t_k]$, to minimize his own cost functional. At $t=t_k$,
the next player (the $(k+1)$-th player, or self--$(k+1)$) takes
over, starting from the initial pair $(t_k,X_k(t_k))$ which is the
{\it terminal pair} of the $k$-th player, and controlling the system
on $[t_k,t_{k+1}]$, etc. Each player knows that the later players
will do their best, and will modify their control systems as well as
their cost functionals. However, in measuring the performance of the
controls, each player will discount the cost/payoff {\it in his/her
own way}. This is the main issue in handling the time-inconsistency,
and it also has to be treated this way so that the results can
recover those for exponential discounting situations. It is expected
that as the mesh size $\|\D\|\equiv\max\{t_k-t_{k-1}\bigm|1\le k\le
N\}\to0$, the Nash equilibrium strategy to the $N$-person
differential game should approach to the desired time-consistent
solution of the original time-inconsistent Problem (N).

\ms

The rest of the paper is organized as follow. In section 2, we
collect some preliminary results, mainly some careful estimates
relevant to our time-inconsistent optimal control problem. Section 3
is devoted to a study of $N$-person differential game. In Section 4,
we will discuss the convergence of Nash equilibrium value function
for the $N$-person differential game, as well as a sufficient
condition for the existence of time-consistent equilibrium control
for Problem (N). Finally, a time-inconsistent LQ problem will be
presented.

\ms

\section{$N$-Person Differential Games}

Consider the following linear controlled ODE parameterized by
$(t,x)\in[0,T]\times\dbR^n$:
\bel{}\left\{\ba{ll}
\ns\ds\dot X(s)=A(t,x,s)X(s)+B(t,x,s)u(s),\qq s\in[t,T],\\
\ns\ds X(t)=x,\ea\right.\ee
with the cost functional
\bel{}\ba{ll}
\ns\ds J(t,x;u(\cd))=\lan G(t,x)X(T),X(T)\ran\\
\ns\ds\qq\qq\qq\q+\int_t^T\big[\lan Q(t,x,s)X(s),X(s)\ran+\lan
R(t,x,s)u(s),u(s)\ran\big]ds.\ea\ee
Here $A$, $B$, $Q$, $R$ and $G$ are some given suitable maps. Let
$\D$ be a partition of $[0,T]$ given by
$$\D:0=t_0<t_1<\cds<t_N=T.$$
We now introduce an $N$-person differential game associated with
$\D$. These $N$ players are labeled by $k=1,2,\cds,N$. The $k$-th
player chooses controls from $\cU[t_{k-1},t_k]$. Any
$(u_1(\cd),\cds,u_N(\cd))\in\cU[t_0,t_1]\times\cds\times\cU[t_{N-1},t_N]$
is identified with $u^\D(\cd)\in\cU[0,T]$ where
\bel{}u^\D(s)=u_k(s),\qq s\in[t_{k-1},t_k),\qq1\le k\le N.\ee
Now, for any $(x,u^\D(\cd))\in\dbR^n\times\cU[0,T]$, let $X^\D(\cd)$
be the solution to the following:
\bel{}\left\{\ba{ll}
\ns\ds\dot
X^\D(s)=A(t_{k-1},X^\D(t_{k-1}),s)X^\D(s)+B(t_{k-1},X^\D(t_{k-1}),s)u^\D(s),\\
\ns\ds\qq\qq\qq\qq\qq s\in(t_{k-1},t_k),\q1\le k\le N,\\
\ns\ds X^\D(0)=x,\ea\right.\ee
The $k$-th player has the following cost functional:
\bel{}\ba{ll}
\ns\ds J_k(u^\D(\cd))\equiv J_k(u_1(\cd),\cds,u_N(\cd))=J(t_{k-1},X^\D(t_{k-1}),u^\D(\cd))\\
\ns\ds\equiv\2n\lan
G(t_{k-1},X^\D(t_{k-1}))X^\D(T),X^\D(T)\ran\\
\ns\ds\q+\int_{t_{k-1}}^T\2n\big[\lan
Q(t_{k-1},X^\D(t_{k-1}),s)X^\D(s),X^\D(s)\ran\1n+\1n\lan
R(t_{k-1},X^\D(t_{k-1}),s)u^\D(s),u^\D(s)\ran\big]ds.\ea\ee

\ms

For any $x\in\dbR^n$ and any partition $\D$ of $[0,T]$, we now pose
the following problem.

\ms

\bf Problem (LQ$^\D$). \rm Find a control $\bar u^\D(\cd)\equiv(\bar
u_1(\cd),\cds,\bar u_N(\cd))\in\cU[0,T]$ such that for each
$k=1,2,\cds,N$,
\bel{}\ba{ll}
\ns\ds J_k(\bar u^\D(\cd))\equiv J_k(\bar u_1(\cd),\cds,\bar
u_{k-1}(\cd),\bar u_k(\cd),\bar
u_{k+1}(\cd),\cds,\bar u_N(\cd))\\
\ns\ds\le J_k(\bar u_1(\cd),\cds,\bar u_{k-1}(\cd),u_k(\cd),\bar
u_{k+1}(\cd),\cds,\bar u_N(\cd)),\qq\forall
u_k(\cd)\in\cU[t_{k-1},t_k],\ea\ee

Any control $\bar u^\D(\cd)$ satisfying the above is called an {\it
equilibrium control} of Problem (LQ$^\D$). The corresponding state
trajectory $\bar X^\D(\cd)$ and the pair $(\bar X^\D(\cd),\bar
u^\D(\cd))$ are called an {\it equilibrium state trajectory} and an
{\it equilibrium pair} of Problem (LQ$^\D$), respectively.

\ms

\ms

We now introduce the following assumptions.

\ms

{\bf(H1)} The maps $A:[0,T]\times[0,T]\to\dbR^{n\times n}$,
$B:[0,T]\times[0,T]\to\dbR^{n\times m}$,
$Q:[0,T]\times[0,T]\to\cS^n$, $R:[0,T]\times[0,T]\to\cS^m$, and
$G:[0,T]\to\cS^n$ are continuous. There exist constants $L,\d>0$
such that
\bel{}\ba{ll}
\ns\ds\|A(t,s)-A(r,s)\|+\|B(t,s)-B(r,s)\|+\|Q(t,s)-Q(r,s)\|\\
\ns\ds+\|R(t,s)-R(r,s)\|+\|G(t)-G(r)\| \le L|t-r|,\qq
s,t,r\in[0,T].\ea\ee
and
\bel{}Q(t,s),G(t)\ge0,\q R(t,s)\ge\d I,\qq\forall t,s\in[0,T].\ee

{\bf(H2)} The maps $G(\cd)$, $Q(\cd\,,\cd)$, and $R(\cd\,,\cd)$
satisfy the following:
\bel{}G(t)\le G(r),\q Q(t,s)\le Q(r,s),\q R(t,s)\le
R(r,s),\qq\forall0\le t\le r\le s\le T.\ee

\ms

For any partition $\D$ of $[0,T]$, we denote
$$\left\{\ba{ll}
\ns\ds A^\D(s)=\sum_{k=1}^NA(t_{k-1},s)I_{[t_{k-1},t_k)}(s),\q
B^\D(s)=\sum_{k=1}^NB(t_{k-1},s)I_{[t_{k-1},t_k)}(s),\\
\ns\ds Q^\D(s)=\sum_{k=1}^NQ(t_{k-1},s)I_{[t_{k-1},t_k)}(s),\q
R^\D(s)=\sum_{k=1}^NR(t_{k-1},s)I_{[t_{k-1},t_k)}(s).\ea\right.$$
Our first result is the following.

\ms

\bf Theorem 2.1. \rm Let (H1) hold. For any partition $\D$ of
$[0,T]$ and any $x\in\dbR^n$, Problem (LQ$^\D$) admits a unique
equilibrium pair $(\bar X^\D(\cd),\bar u^\D(\cd))$. Moreover, $\bar
X^\D(\cd)$ and $\bar u^\D(\cd)$ are linked by the following:
\bel{feedback}\bar u^\D(s)=-R^\D(s)^{-1}B^\D(s)^TP^\D(s)\bar
X^\D(s),\qq s\in[0,T],\ee
where $P^\D(\cd)$ is the unique solution to the following Riccati
equation:
\bel{Riccati}\left\{\ba{ll}
\ns\ds\dot
P^\D(s)+P^\D(s)A^\D(s)+A^\D(s)^TP^\D(s)+Q^\D(s)\\
\ns\ds\qq\qq-P^\D(s)B^\D(s)R^\D(s)^{-1}B^\D(s)^TP^\D(s)=0,\qq
s\in(t_{k-1},t_k),\\
\ns\ds
P^\D(t_k-0)=\F^\D(t_N;t_k)^TG(t_{k-1})\F^\D(t_N;t_k)\\
\ns\ds\qq\qq+\int_{t_k}^{t_N}\(\F^\D(s;t_k)^TQ(t_{k-1},s)\F^\D(s;t_k)
+\Psi^\D(s;t_k)^TR(t_{k-1},s)\Psi^\D(s;t_k)\)ds,\\
\ns\ds\qq\qq\qq\qq\qq\qq\qq\qq\qq\qq1\le k\le N,\ea\right.\ee
with $\F^\D(\cd\,;t_k)$ ($0\le k\le N-1$) being the solution to the
following:
\bel{F}\left\{\ba{ll}
\ns\ds\F^\D_s(s;t_k)=\[A^\D(s)-B^\D(s)R^\D(s)^{-1}B^\D(s)^TP^\D(s)\]
\F^\D(s;t_k),\q
s\in(t_k,T],\\
\ns\ds\F^\D(t_k;t_k)=I,\ea\right.\ee
and
\bel{Psi}\Psi^\D(s;t_k)=-R^\D(s)^{-1}B^\D(s)P^\D(s)\F^\D(s;t_k),\qq
s\in[t_k,T].\ee
The equilibrium state trajectory $\bar X^\D(\cd)$ is the solution to
the following closed-loop system:
\bel{}\left\{\ba{ll}
\ns\ds\dot{\bar
X}^\D(s)=\[A^\D(s)-B^\D(s)R^\D(s)^{-1}B^\D(s)^TP^\D(s)\]\bar
X^\D(s),\qq s\in[0,T],\\
\ns\ds\bar X(0)=x,\ea\right.\ee
and the equilibrium pair $(\bar X^\D(\cd),\bar u^\D(\cd))$ can be
explicitly represented by the following:
\bel{}\left\{\ba{ll}
\ns\ds\bar X^\D(s)=\F^\D(s;0)x,\\
\ns\ds\bar u^\D(s)=\Psi^\D(s;0)x,\ea\right.\qq s\in[0,T].\ee
Moreover,
\bel{}0\le P^\D(t)\le P^\D_0(t),\qq t\in[0,T],\ee
where $P^\D_0(\cd)$ is the unique solution to the following Lyapunov
equation:
\bel{Lyapunov}\left\{\ba{ll}
\ns\ds\dot
P^\D_0(s)+P^\D_0(s)A^\D(s)+A^\D(s)^TP^\D_0(s)+Q^\D(s)=0,\qq
s\in(t_{k-1},t_k),\\
\ns\ds
P^\D_0(t_k-0)=\F^\D(t_N;t_k)^TG(t_{k-1})\F^\D(t_N;t_k)\\
\ns\ds\qq\qq+\int_{t_k}^{t_N}\big[\F^\D(s;t_k)^TQ(t_{k-1},s)\F^\D(s;t_k)
+\Psi^\D(s;t_k)^TR(t_{k-1},s)\Psi^\D(s;t_k)\big]ds,\\
\ns\ds\qq\qq\qq\qq\qq\qq\qq\qq\qq\qq1\le k\le N,\ea\right.\ee

\ms

We point out that the solution $P^\D(\cd)$ of Riccati equation
(\ref{Riccati}) and the soluion$P^\D_0(\cd)$ of Lyapunov equation
(\ref{Lyapunov}) have possible jumps at $t_k$, $k=1,2,\cds,N-1$.

\ms

\it Proof. \rm Let $x\in\dbR^n$ and $\D:0=t_0<t_1<\cds<t_N=T$ be
given. Let $(\bar X^\D(\cd),\bar u^\D(\cd))$ be an equilibrium pair
of Problem (LQ$^\D$). Then the restriction of which on
$[t_{N-1},t_N]$ is the optimal pair of the LQ problem for Player $N$
on $[t_{N-1},t_N]$, with the state equation
\bel{}\left\{\ba{ll}
\ns\ds\dot X^\D(s)=A^\D(s)X^\D(s)+B^\D(s)u_N(s),\qq s\in[t_{N-1},t_N],\\
\ns\ds X^\D(t_{N-1})=\bar X^\D(t_{N-1}),\ea\right.\ee
and with the cost functional
\bel{}\ba{ll}
\ns\ds J_N(\bar u_1(\cd),\cds,\bar u_{N-1}(\cd),u_N(\cd))=\lan
G_NX^\D(t_N),X^\D(t_N)\ran\\
\ns\ds\qq\qq\qq\qq\qq\qq+\int_{t_{N-1}}^{t_N}\big[\lan
Q^\D(s)X^\D(s),X^\D(s)\ran+\lan
R^\D(s)u_N(s),u_N(s)\ran\big]ds,\ea\ee
where $G_N=G(t_{N-1})$. To study this LQ problem, we consider the
following state equation:
\bel{}\left\{\ba{ll}
\ns\ds\dot X_N(s)=A^\D(s)X_N(s)+B^\D(s)u_N(s),\qq s\in[t,t_N],\\
\ns\ds X_N(t)=y,\ea\right.\ee
where $(t,y)\in[t_{N-1},t_N)\times\dbR^n$, with the cost functional
\bel{}\ba{ll}
\ns\ds J_N(t,y;u_N(\cd))=\lan
G_NX_N(t_N),X_N(t_N)\ran\\
\ns\ds\qq\qq\qq\qq+\int_t^{t_N}\big[\lan
Q^\D(s)X_N(s),X_N(s)\ran+\lan R^\D(s)u_N(s),u_N(s)\ran\big]ds,\ea\ee
For such an LQ problem on $[t,t_N]$, under (H1), there exists a
unique optimal control which must have the following form:
\bel{u_N}\bar u_N(s)=-R^\D(s)^{-1}B^\D(s)^TP^\D(s)\bar X_N(s),\qq
s\in[t,t_N],\ee
where $P^\D(\cd)$ is the unique solution to the following Riccati
equation:
\bel{3.19}\left\{\ba{ll}
\ns\ds\dot
P^\D(s)+P^\D(s)A^\D(s)+A^\D(s)^TP^\D(s)+Q^\D(s)\\
\ns\ds\qq\qq-P^\D(s)B^\D(s)R^\D(s)^{-1}B^\D(s)^TP^\D(s)=0,\qq
s\in(t_{N-1},t_N),\\
\ns\ds P^\D(t_N)=G_N,\ea\right.\ee
and $\bar X_N(\cd)$ is the solution to the following closed-loop
state equation:
\bel{}\left\{\ba{ll}
\ns\ds\dot{\bar
X}_N(s)=\[A^\D(s)-B^\D(s)R^\D(s)^{-1}B^\D(s)^TP^\D(s)\]\bar
X_N(s),\qq s\in[t,t_N],\\
\ns\ds\bar X_N(t)=y.\ea\right.\ee
Let $\F^\D(\cd\,;t)$ be the solution to the following: (note that
$t\in[t_{N-1},t_N]$)
\bel{2.26}\left\{\ba{ll}
\ns\ds\F^\D_s(s;t)=\[A^\D(s)-B^\D(s)R^\D(s)^{-1}B^\D(s)^TP^\D(s)\]
\F^\D(s;t),\q
s\in(t,t_N],\\
\ns\ds\F^\D(t;t)=I,\ea\right.\ee
and denote
\bel{}\Psi^\D(s;t)=-R^\D(s)^{-1}B^\D(s)P^\D(s)\F^\D(s;t),\qq
s\in[t,t_N].\ee
Then the optimal pair $(\bar X_N(\cd),\bar u_N(\cd))$ of LQ problem
(on $[t,t_N]$) admits the following representation:
\bel{}\left\{\ba{ll}
\ns\ds\bar X_N(s)=\F^\D(s;t)y,\\
\ns\ds\bar u_N(s)=\Psi^\D(s;t)y,\ea\qq s\in[t,t_N].\right.\ee
Further,
\bel{}\ba{ll}
\ns\ds\lan P^\D(t)y,y\ran=J_N(t,y;\bar u_N(\cd))=\lan G_N\bar
X_N(t_N),\bar
X_N(t_N)\ran\\
\ns\ds\qq\qq\qq+\int_t^{t_N}\big[\lan Q(t_{N-1},s)\bar X_N(s),\bar
X_N(s)\ran+\lan R(t_{N-1},s)\bar u_N(s),\bar
u_N(s)\ran\big]ds\\
\ns\ds=\lan\big[\F^\D(t_N;t)^TG_N\F^\D(t_N;t)+\int_t^{t_N}\big(\F^\D(s;t)^TQ(t_{N-1},s)
\F^\D(s;t)\\
\ns\ds\qq\qq\qq+\Psi^\D(s;t)^TR(t_{N-1},s)\Psi^\D(s;t)\big)ds\big]y,y\ran.\ea\ee
Since $y\in\dbR^n$ can be arbitrarily chosen, we have
\bel{}\ba{ll}
\ns\ds
P^\D(t)=\F^\D(t_N;t)^TG_N\F^\D(t_N;t)+\int_t^{t_N}\big[\F^\D(s;t)^TQ(t_{N-1},s)
\F^\D(s;t)\\
\ns\ds\qq\qq\qq\qq\qq\qq+\Psi^\D(s;t)^TR(t_{N-1},s)\Psi^\D(s;t)\big]ds,\q
t\in(t_{N-1},t_N].\ea\ee
Also, by the optimality of $\bar u_N(\cd)$, we have
\bel{}\ba{ll}
\ns\ds\lan P^\D(t)y,y\ran=J_N(t,y;\bar u_N(\cd))\le J_N(t,y;0)\\
\ns\ds\qq=\lan G_NX^0(t_N),X^0(t_N)\ran+\int_t^{t_N}\lan
Q(t_{N-1},s)X^0(s),X^0(s)\ran ds=\lan P^\D_0(t)y,y\ran,\ea\ee
where $X^0(\cd)$ is the solution to the following:
\bel{}\left\{\ba{ll}
\ns\ds\dot X^0(s)=A^\D(s)X^0(s),\qq s\in[t,t_N],\\
\ns\ds X^0(t)=y,\ea\right.\ee
and $P^\D_0(\cd)$ is the solution to the following Lyapunov
equation:
\bel{}\left\{\ba{ll}
\ns\ds\dot
P^\D_0(s)+P^\D_0(s)A^\D(s)+A^\D(s)^TP^\D_0(s)+Q^\D(s)=0,\qq
s\in(t_{N-1},t_N),\\
\ns\ds P^\D_0(t_N-0)=G_N,\ea\right.\ee
which can be represented by the following:
\bel{}\ba{ll}
\ns\ds
P^\D_0(t)=\F^\D_0(t_N;t)^TG_N\F^\D_0(t_N;t)+\int_t^{t_N}\F^\D_0(s,t)^TQ(t_{N-1},s)
\F^\D_0(s;t)ds,\q t\in[t,t_N],\ea\ee
with $\F^\D_0(\cd\,;t)$ being the solution to the following:
\bel{}\left\{\ba{ll}
\ns\ds{\pa\over\pa s}\F^\D_0(s;t)=A^\D(s)\F^\D_0(s;t),\qq
s\in[t,t_N],\\
\ns\ds\F^\D_0(t;t)=I.\ea\right.\ee
Note that $\F^\D_0(\cd\,;t)$ can be defined for any $t\in[0,t_N)$,
which will be used below. Hence,
\bel{}0\le P^\D(t)\le P^\D_0(t),\qq t\in(t_{N-1},t_N].\ee
It is also clear that the restriction of the equilibrium pair $(\bar
X^\D(\cd),\bar u^\D(\cd))$ on $(t_{N-1},t_N]$ admits the following
representation:
\bel{}\left\{\ba{ll}
\ns\ds\bar X^\D(s)=\F^\D(s;t_{N-1})\bar X^\D(t_{N-1}),\\
\ns\ds\bar u^\D(s)=\Psi^\D(s;t_{N-1})\bar X^\D(t_{N-1}),\ea\qq
s\in[t_{N-1},t_N].\right.\ee

\ms

Next, for Player $(N-1)$, inspired by the above, we consider the
following state equation:
\bel{2.38}\left\{\ba{ll}
\ns\ds\dot X_{N-1}(s)=A^\D(s)X_{N-1}(s)+B^\D(s)u_{N-1}(s),\qq s\in[t,t_{N-1}],\\
\ns\ds X_{N-1}(t)=y,\ea\right.\ee
where $(t,y)\in[t_{N-2},t_{N-1})\times\dbR^n$. Let us denote
\bel{}\left\{\ba{ll}
\ns\ds\wt X^\D_{N-1}(s)=\F^\D(s;t_{N-1})X_{N-1}(t_{N-1}),\\
\ns\ds\wt u^\D_{N-1}(s)=\Psi^\D(s;t_{N-1})X_{N-1}(t_{N-1}),\ea\qq
s\in[t_{N-1},t_N].\right.\ee
Thus, $(\wt X^\D_{N-1}(\cd),\wt u^\D_{N-1}(\cd))$ is the optimal
pair for Player $N$ starting from the initial pair
$(t_{N-1},X_{N-1}$ $(t_{N-1}))$. The cost functional for the LQ
problem of Player $(N-1)$ on $[t,t_{N-1}]$ is taken to be
\bel{2.40}\ba{ll}
\ns\ds J_{N-1}(t,y;u_{N-1}(\cd))\\
\ns\ds=\int_t^{t_{N-1}}\big[\lan
Q(t_{N-2},s)X_{N-1}(s),X_{N-1}(s)\ran+\lan
R(t_{N-2},s)u_{N-1}(s),u_{N-1}(s)\ran\big]ds\\
\ns\ds\qq+\int_{t_{N-1}}^{t_N}\big[\lan Q(t_{N-2},s)\wt
X^\D_{N-1}(s),\wt X^\D_{N-1}(s)\ran+\lan
R(t_{N-2},s)\wt u^\D_{N-1}(s),\wt u^\D_{N-1}(s)\ran\big]ds\\
\ns\ds\qq+\lan G(t_{N-2})\wt X^\D(t_N),\wt X^\D(t_N)\ran\\
\ns\ds\equiv\int_t^{t_{N-1}}\big[\lan
Q^\D(s)X_{N-1}(s),X_{N-1}(s)\ran+\lan
R^\D(s)u_{N-1}(s),u_{N-1}(s)\ran\big]ds\\
\ns\ds\qq+\lan G_{N-1}X_{N-1}(t_{N-1}),X_{N-1}(t_{N-1})\ran,\ea\ee
where
\bel{}\ba{ll}
\ns\ds
G_{N-1}=\F^\D(t_N;t_{N-1})^TG(t_{N-2})\F^\D(t_N;t_{N-1})\\
\ns\ds\qq+\2n\int_{t_{N-1}}^{t_N}\3n\2n\big[\F^\D(s;t_{N-1})^T\1n
Q(t_{N-2},s)\F^\D(s;t_{N-1})\2n+\2n\Psi^\D(s;t_{N-1})^T\1n
R(t_{N-2},s)\Psi^\D(s;t_{N-1})\big]ds.\ea\ee
For such an LQ problem (on $[t,t_{N-1}]$), under (H1), the optimal
control is given by
\bel{}\bar u_{N-1}(s)=-R^\D(s)^{-1}B^\D(s)^TP^\D(s)\bar
X_{N-1}(s),\qq s\in[t,t_{N-1}],\ee
where $P^\D(\cd)$ is the solution to the following Riccati equation:
\bel{}\left\{\ba{ll}
\ns\ds\dot
P^\D(s)+P^\D(s)A^\D(s)+A^\D(s)^TP^\D(s)+Q^\D(s)\\
\ns\ds\qq\qq-P^\D(s)B^\D(s)R^\D(s)^{-1}B^\D(s)^TP^\D(s)=0,\qq
s\in(t_{N-2},t_{N-1}),\\
\ns\ds P^\D(t_{N-1}-0)=G_{N-1},\ea\right.\ee
and $\bar X_{N-1}(\cd)$ is the solution to the following closed-loop
state equation:
\bel{}\left\{\ba{ll}
\ns\ds\dot{\bar
X}_{N-1}(s)=\[A^\D(s)-B^\D(s)R^\D(s)^{-1}B^\D(s)^TP^\D(s)\]\bar
X_{N-1}(s),\qq s\in[t,t_{N-1}],\\
\ns\ds\bar X_{N-1}(t)=y.\ea\right.\ee
Now, similar to (\ref{2.26}), for $t\in[t_{N-2},t_{N-1}]$, let
$\F^\D(\cd\,;t)$ be the solution to the following:
\bel{}\left\{\ba{ll}
\ns\ds\F^\D_s(s;t)=\[A^\D(s)-B^\D(s)R^\D(s)^{-1}B^\D(s)^TP^\D(s)\]
\F^\D(s;t),\q
s\in(t,t_N],\\
\ns\ds\F^\D(t;t)=I,\ea\right.\ee
and denote
\bel{}\Psi^\D(s;t)=-R^\D(s)^{-1}B^\D(s)P^\D(s)\F^\D(s;t),\qq
s\in[t,t_N].\ee
It is clear that for $t\in[t_{N-2},t_{N-1}]$,
\bel{}\left\{\ba{ll}
\ns\ds\F^\D(s;t)=\F^\D(s;t_{N-1})\F^\D(t_{N-1};t),\\
\ns\ds\Psi^\D(s;t)=\Psi^\D(s;t_{N-1})\F^\D(t_{N-1};t),\ea\right.\qq
s\in[t_{N-1},t_N],\ee
and the optimal pair $(\bar X_{N-1}(\cd),\bar u_{N-1}(\cd))$ of the
LQ problem associated with (\ref{2.38}) and (\ref{2.40}) (on
$[t,t_{N-1}]$) is given by the following:
\bel{}\left\{\ba{ll}
\ns\ds\bar X_{N-1}(s)=\F^\D(s;t)y,\\
\ns\ds\bar u_{N-1}(s)=\Psi^\D(s;t)y,\ea\right.\qq
s\in[t,t_{N-1}].\ee
Hence, the restriction of the equilibrium pair $(\bar X^\D(\cd),\bar
u^\D(\cd))$ on $[t_{N-2},t_N]$ admits the following representation:
\bel{}\left\{\ba{ll}
\ns\ds\bar X^\D(s)=\F^\D(s;t_{N-2})\bar X^\D(t_{N-2}),\\
\ns\ds\bar u^\D(s)=\Psi^\D(s;t_{N-2})\bar
X^\D(t_{N-2}),\ea\right.\qq s\in[t_{N-2},t_N].\ee
Further,
\bel{}\ba{ll}
\ns\ds\lan P^\D(t)y,y\ran=J_{N-1}(t,y;\bar u_{N-1}(\cd))=\lan
G_{N-1}\bar X_{N-1}(t_{N-1}),\bar
X_{N-1}(t_{N-1})\ran\\
\ns\ds\qq\qq+\int_t^{t_{N-1}}\big[\lan Q(t_{N-2},s)\bar
X_{N-1}(s),\bar X_{N-1}(s)\ran+\lan R(t_{N-2},s)\bar u_{N-1}(s),\bar
u_{N-1}(s)\ran\big]ds\\
\ns\ds=\lan\big[\F^\D(t_{N-1};t)^TG_{N-1}\F^\D(t_{N-1};t)+\int_t^{t_{N-1}}
\big[\F^\D(s;t)^TQ(t_{N-2},s)\F^\D(s;t)\\
\ns\ds\qq\qq\qq+\Psi^\D(s;t)^TR(t_{N-2},s)\Psi^\D(s;t)\big]ds\big]y,y\ran\\
\ns\ds=\lan\big[\F^\D(t_N;t)^TG(t_{N-2})\F^\D(t_N;t)+\int_t^{t_N}\big(\F^\D(s;t)^TQ(t_{N-2},s)
\F^\D(s;t)\\
\ns\ds\qq\qq\qq+\Psi^\D(s;t)^TR(t_{N-2},s)\Psi^\D(s;t)\big)ds\big]y,y\ran.\ea\ee
Since $y\in\dbR^n$ can be arbitrarily chosen, we have
\bel{}\ba{ll}
\ns\ds
P^\D(t)=\F^\D(t_N;t)^TG(t_{N-2})\F^\D(t_N;t)+\int_t^{t_N}\big[\F^\D(s;t)^TQ(t_{N-2},s)
\F^\D(s;t)\\
\ns\ds\qq\qq\qq\qq\qq\qq+\Psi^\D(s;t)^TR(t_{N-2},s)\Psi^\D(s;t)\big]ds,\q
t\in(t_{N-2},t_{N-1}).\ea\ee
Also, by the optimality of $\bar u_{N-1}(\cd)$, we have
\bel{}\ba{ll}
\ns\ds\lan P^\D(t)y,y\ran=J_{N-1}(t,y;\bar u_{N-1}(\cd))\le J_{N-1}(t,y;0)\\
\ns\ds\q=\lan
G_{N-1}X^0(t_{N-1}),X^0(t_{N-1})\ran+\int_t^{t_{N-1}}\lan
Q(t_{N-2},s)X^0(s),X^0(s)\ran ds=\lan P^\D_0(t)y,y\ran,\ea\ee
where $X^0(\cd)$ is the solution to the following:
\bel{}\left\{\ba{ll}
\ns\ds\dot X^0(s)=A^\D(s)X^0(s),\qq s\in[t,t_{N-1}],\\
\ns\ds X^0(t)=y,\ea\right.\ee
and $P^\D_0(\cd)$ is the solution to the following Lyapunov
equation:
\bel{}\left\{\ba{ll}
\ns\ds\dot
P^\D_0(s)+P^\D_0(s)A^\D(s)+A^\D(s)^TP^\D_0(s)+Q^\D(s)=0,\qq
s\in(t_{N-2},t_{N-1}),\\
\ns\ds P^\D_0(t_{N-1}-0)=G_{N-1},\ea\right.\ee
which, similar to the above, admits the following representation:
\bel{}\ba{ll}
\ns\ds
P^\D_0(t)\1n=\1n\F^\D_0(t_{N-1};t)^TG_{N-1}\F^\D_0(t_{N-1};t)\1n+\2n\int_t^{t_{N-1}}\2n\F^\D_0(s,t)^TQ^\D(s)
\F^\D_0(s;t)ds,\q t\in[t,t_{N-1}],\ea\ee
Hence,
\bel{}0\le P^\D(t)\le P^\D_0(t),\qq t\in(t_{N-2},t_{N-1}).\ee
Then one can apply induction to complete the proof. \endpf

\ms

\section{Time-Consistent Solutions}

We now pose the following problem.

\ms

\bf Problem (LQ). \rm For any given $x\in\dbR^n$, find a control
$\bar u(\cd)\in\cU[0,T]$ satisfying the following: For any $\e>0$,
there exists a $\d>0$ such that for any partition $\D$ of $[0,T]$
with $\|\D\|<\d$, one has
\bel{}J_k(\bar u(\cd))\le J_k(\bar u^\D(\cd))+\e.\ee

\ms

Any control $\bar u(\cd)\in\cU[0,T]$ satisfying the above is called
an {\it equilibrium control} of Problem (LQ). The corresponding
state trajectory $\bar X(\cd)$ and the pair $(\bar X(\cd),\bar
u(\cd))$ are called an {\it equilibrium state trajectory} and an
{\it equilibrium pair} of Problem (LQ), respectively.

\ms

The following gives a weaker notion of time-consistent solutions to
Problem (LQ).

\ms

\bf Definition 3.1. \rm A control $\bar u(\cd)\in\cU[0,T]$ is called
a {\it weak equilibrium control} of Problem (LQ) if there exists a
sequence of partitions $\D_m$ of $[0,T]$ with $\|\D_m\|\to0$ so that
for any $\e>0$, there exists an $m_0>0$ such that
\bel{}J_k(\bar u(\cd))\le J_k(\bar u^{\D_m}(\cd))+\e,\qq\forall m\ge
m_0.\ee

Our next goal is to find the limit as the mesh size $\|\D\|$ of $\D$
approaches to zero. For this, we need (H2).

\ms

\bf Theorem 3.1. \rm Let (H1)--(H2) hold. Then for any partition
$\D$ of $[0,T]$,
\bel{}0\le P^\D(t)\le P^\D_0(t)\le\bar P^\D_0(t),\qq t\in[0,T],\ee
where $\bar P^\D_0(\cd)$ is the unique solution to the following
Lyapunov equation:
\bel{Lyapunov2}\left\{\ba{ll}
\ns\ds\dot{\bar P}_0^\D(s)+\bar P^\D_0(s)A^\D(s)+A^\D(s)^T\bar
P^\D_0(s)+Q^\D(s)=0,\qq
s\in(0,t_N),\\
\ns\ds \bar P^\D_0(t_N)=G(t_{N-1}).\ea\right.\ee
Consequently, $P^\D(\cd)$ is bounded uniformly in $\D$.

\ms

\it Proof. \rm Recall that for $k=1,2,\cds,N-1$,
\bel{}\ba{ll}\ns\ds
P^\D_0(t_k-0)=\F^\D(t_N;t_k)^TG(t_{k-1})\F^\D(t_N;t_k)\\
\ns\ds\qq\qq+\int_{t_k}^{t_N}\big[\F^\D(s;t_k)^TQ(t_{k-1},s)\F^\D(s;t_k)
+\Psi^\D(s;t_k)^TR(t_{k-1},s)\Psi^\D(s;t_k)\big]ds,\ea\ee
and (making use of the monotonicity of $t\mapsto G(t)$, $t\mapsto
Q(s,t)$, and $t\mapsto R(s,t)$)
\bel{2.60}\ba{ll}
\ns\ds P^\D_0(t_k+0)\ge
P^\D(t_k+0)=\F^\D(t_N;t_k)^TG(t_k)\F^\D(t_N;t_k)\\
\ns\ds\qq\qq+\int_{t_k}^{t_N}\big[\F^\D(s;t_k)^T
Q(t_k,s)\F^\D(s;t_k)+\Psi^\D(s;t_k)^TR(t_k,s)\Psi^\D(s;t_k)\big]ds\\
\ns\ds\ge\F^\D(t_N;t_k)^TG(t_{k-1})\F^\D(t_N;t_k)+\int_{t_k}^{t_N}\big[\F^\D(s;t_k)^T
Q(t_{k-1},s)\F^\D(s;t_k)\\
\ns\ds\qq\qq\qq\qq+\Psi^\D(s;t_k)^TR(t_{k-1},s)\Psi^\D(s;t_k)\big]ds=P^\D(t_k-0)=P^\D_0(t_k-0).\ea\ee
Note that
\bel{}P^\D_0(t)=\bar P^\D_0(t),\qq t\in(t_{N-1},t_N],\ee
and due to (making use of (\ref{2.60}) for $k=N-1$)
\bel{}P^\D_0(t_{N-1}-0)\le P^\D(t_{N-1}+0)\le P^\D_0(t_{N-1}+0)=\bar
P^\D_0(t_{N-1})=\bar P^\D_0(t_{N-1}-0),\ee
we have
\bel{}P^\D_0(t)\le\bar P^\D_0(t),\qq t\in(t_{N-2},t_{N-1}).\ee
Then by induction, we can obtain
\bel{}P^\D_0(t)\le\bar P^\D_0(t),\qq t\in[0,t_N].\ee
By the boundness of $A(\cd\,,\cd)$ and $Q(\cd\,,\cd)$, we have the
boundness of $\bar P^\D_0(\cd)$ uniformly in $\D$. Hence, we
complete the proof. \endpf

\ms

We see that $P^\D(\cd)$ has a possible jump at each $t_k$, with the
jump size:
\bel{}\ba{ll}
\ns\ds\D P^\D(t_k)\equiv P^\D(t_k+0)-P^\D(t_k-0)\\
\ns\ds\qq\q~~~=\F^\D(t_N;t_k)^T\big[G(t_k)-G(t_{k-1})\big]\F^\D(t_N;t_k)\\
\ns\ds\qq\qq\qq+\int_{t_k}^{t_N}\big[\F^\D(s;t_k)^T
\big(Q(t_k,s)-Q(t_{k-1},s)\big)\F^\D(s;t_k)\\
\ns\ds\qq\qq\qq\qq+\Psi^\D(s;t_k)^T\big(R(t_k,s)-R(t_{k-1},s)\big)
\Psi^\D(s;t_k)\big]ds\ge0.\ea\ee
By (H1)--(H2), we have
\bel{}\|\D P^\D(t_k)\|\le K(t_k-t_{k-1})\le K\|\D\|.\ee
Next, we define $\wt P^\D(\cd)$ as follows:
\bel{}\left\{\ba{ll}
\ns\ds\wt P^\D(t)=P^\D(t),\qq t\in(t_{N-1},t_N],\\
\ns\ds\wt P^\D(t)=P^\D(t)+{t-t_{k-1}\over t_k-t_{k-1}}\D
P^\D(t_k),\qq t\in(t_{k-1},t_k),\\
\ns\ds\wt P^\D(t_k)=P^\D(t_k+0),\qq1\le k\le N-1.\ea\right.\ee
Then $\{\wt P^\D(\cd)\}$ is uniformly bounded and equicontinuous.
Hence, we may assume that along a certain sequence $\D_m$ with
$\|\D_m\|\to0$,
\bel{}\lim_{m\to\infty}\wt P^{\D_m}(\cd)=P(\cd),\ee
for some $P(\cd)$. Also, we have
\bel{}\|\wt P^\D(t)-P^\D(t)\|\le\max_{1\le k\le N-1}\|\D
P^\D(t_k)\|\le K\|\D\|\to0,\qq\hb{as }\|\D\|\to0.\ee
Hence, we have
\bel{}\lim_{m\to\infty}\|P^{\D_m}(\cd)-P(\cd)\|=0.\ee
Next, it is clear that
\bel{}\lim_{\|\D\|\to0}\2n\Big\{\1n\|A^\D(s)-A(s,s)\|\1n+\1n\|B^\D(s)-B(s,s)\|\1n+\1n
\|Q^\D(s)-Q(s,s)\|\1n+\1n\|R^\D(s)-R(s,s)\|\1n\Big\}\2n=0.\ee
Hence,
\bel{}\lim_{\|\D\|\to0}\|\F^\D(s;t)-\F(s;t)\|=0,\ee
with $\F(\cd\,;t)$ being the solution to the following:
\bel{}\left\{\ba{ll}
\ns\ds\F_s(s;t)=\[A(s,s)-B(s,s)R(s,s)^{-1}B(s,s)^TP(s)\]\F(s;t),\q
s\in(t,t_N],\\
\ns\ds\F(t;t)=I.\ea\right.\ee
Consequently, $P(\cd)$ satisfies the following:
\bel{}\ba{ll}
\ns\ds P(t)=\F(T;t)^T\1n G(t)\F(T;t)\1n+\1n\int_t^T\big[\F(s;t)^T\1n
Q(t,s)\F(s;t)\\
\ns\ds\qq\q+\F(s;t)^T\1n P(s)B(s,s)^T\1n R(s,s)^{-1}\1n
R(t,s)R(s,s)^{-1}\1n B(s,s)P(s)\F(s;t)\big]ds,\q t\in(0,T).\ea\ee
Denote
$$A(s)=A(s,s),\q B(s)=B(s,s),\q R(s)=R(s,s).$$
Then we have the following system of forward-backward Volterra
integral equations:
\bel{}\left\{\ba{ll}
\ns\ds\F(s;t)=I+\int_t^s\big[A(r)-B(r)R(r)^{-1}B(r)^TP(r)\big]\F(r;t)dr,\qq s\in[t,T],\\
\ns\ds P(t)=\F(T;t)^T\1n G(t)\F(T;t)\1n+\1n\int_t^T\big[\F(r;t)^T\1n
Q(t,r)\F(r;t)\\
\ns\ds\qq\q+\F(r;t)^T\1n P(r)B(r)^T\1n R(r)^{-1}\1n
R(t,r)R(r)^{-1}\1n B(r)P(r)\F(r;t)\big]dr,\q t\in[0,T].\ea\right.\ee

Suppose the above admits a unique solution $(\F(\cd\,;\cd),P(\cd))$.
Then
\bel{}\lim_{\|\D\|\to0}\|P^\D(\cd)-P(\cd)\|=0,\ee
and
\bel{}\lim_{\|\D\|\to0}\|\F^\D(\cd\,;\cd)-\F(\cd\,;\cd)\|=0.\ee
Then
\bel{}\lim_{\|\D\|\to0}\Big\{\|\bar X^\D(\cd)-\bar X(\cd)\|+\|\bar
u^\D(\cd)-\bar u(\cd)\|\Big\}=0,\ee
with
\bel{}\left\{\ba{ll}
\ns\ds\bar X(s)=\F(s;0)x,\\
\ns\ds\bar u(s)=-R(s)^{-1}B(s)P(s)\bar
X(s)\equiv-R(s)^{-1}B(s)P(s)\F(s;0)x,\ea\right.\qq s\in[0,T].\ee
Hence, for any $\e>0$, there exists a $\d>0$ such that for any
partition $\D$ of $[0,T]$, as long as $\|\D\|<\d$, one has
\bel{}J_k(\bar u(\cd))\le J_k(\bar u^\D(\cd))+K\|\D\|<J_k(\bar
u^\D(\cd))+\e.\ee

 \ms

\bs

\no{\Large\bf Appendix.}

\ms

\rm

Let us now solve Problem (C) explicitly. For any given
$(t,x)\in[0,T)\times\dbR^n$, according to a standard LQ theory, in
the current case, the corresponding Riccati differential equation
reads
$$\left\{\ba{ll}
\ns\ds\dot P(s)-P(s)^2=0,\qq s\in[t,T],\\
\ns\ds P(T)=h(t),\ea\right.
\eqno(A.1)$$
Clearly the solution of the above Riccati equation depends on $t$.
Hence, we denote it by $P(\cd\,;t)$. Simple calculation shows that
$$P(s;t)={h(t)\over1+h(t)(T-s)},\qq s\in[t,T].
\eqno(A.2)$$
The optimal control trajectory is the solution to the following
closed-loop system
$$\left\{\ba{ll}
\ns\ds\dot{\bar X}(s)=-P(s;t)\bar X(s),\qq s\in[t,T],\\
\ns\ds\bar X(t)=x,\ea\right.
\eqno(A.3)$$
which is given by
$$\bar X(s;t,x)=x{1+h(t)(T-s)\over1+h(t)(T-t)}\,,\qq
s\in[t,T],
\eqno(A.4)$$
and the optimal control is given by
$$\bar u(s;t,x)=-P(s;t)\bar
X(s;t,x)=-{xh(t)\over1+h(t)(T-t)},\qq s\in[t,T].
\eqno(A.5)$$
Now, if we let
$$J(t;\t,y;u(\cd))=\int_\t^Tu(s)^2ds+h(t)X(T;\t,y,u(\cd))^2,\qq
\t\in[t,T],
\eqno(A.6)$$
then the optimal value function (for fixed $t$) is given by
$$\ba{ll}
\ns\ds V(t;\t,y)\equiv\inf_{u(\cd)\in\cU[\t,T]}
J(t;\t,y;u(\cd))=J(t;\t,y;\bar u(\cd))=P(\t;t)y^2\\
\ns\ds\qq\q~~={h(t)\over
1+h(t)(T-\t)}y^2,\qq\forall(\t,y)\in[t,T]\times\dbR.\ea
\eqno(A.7)$$
Next, let $\t\in(t,T)$, we consider Problem (C) on $[\t,T]$ with
initial state
$$y=\bar X(\t;t,x)=x{1+h(t)(T-\t)\over1+h(t)(T-t)}.
\eqno(A.8)$$
The same as above, we see that the corresponding solution to the
Riccati equation is given by
$$P(s;\t)={h(\t)\over1+h(\t)(T-s)},\qq s\in[\t,T],
\eqno(A.9)$$
and
$$\ba{ll}
\ns\ds\inf_{u(\cd)\in\cU[\t,T]}J(\t;\t,y;u(\cd))=P(\t;\t)y^2={h(\t)y^2\over
1+h(\t)(T-\t)}\\
\ns\ds={x^2h(\t)[1+h(t)(T-\t)]^2\over[1+h(\t)(T-\t)][1+h(t)(T-t)]^2}.\ea
\eqno(A.10)$$
However,
$$\ba{ll}
\ns\ds J(\t,y;\bar u(\cd))=\int_\t^T\bar u(s)^2ds+h(\t)
X(T;\t,y,\bar u(\cd))^2\\
\ns\ds={x^2h(t)^2(T-\t)\over[1+h(t)(T-t)]^2}+h(\t)\[y-{xh(t)(T-\t)\over1+h(t)(T-t)}\]^2\\
\ns\ds={x^2h(t)^2(T-\t)\over[1+h(t)(T-t)]^2}+{x^2h(\t)\over[1+h(t)(T-t)]^2}.\ea
\eqno(A.11)$$
Hence,
$$\ba{ll}
\ns\ds J(\t,y;\bar
u(\cd))-\inf_{u(\cd)\in\cU[\t,T]}J(\t,y;u(\cd))\\
\ns\ds={x^2h(t)^2(T-\t)+x^2h(\t)\over[1+h(t)(T-t)]^2}-{x^2h(\t)[1+h(t)(T-\t)]^2\over[1+h(\t)(T-\t)]
[1+h(t)(T-t)]^2}\\
\ns\ds=x^2\[{[h(t)^2(T-\t)+h(\t)][1+h(\t)(T-\t)]-h(\t)[1+h(t)(T-\t)]^2\over[1+h(\t)(T-\t)][1+h(t)(T-t)]^2}\]\\
\ns\ds={x^2[h(\t)-h(t)]^2(T-\t)\over[1+h(\t)(T-\t)][1+h(t)(T-t)]^2}>0,\qq\qq\hb{unless
$h(\t)=h(t)$ or $x=0$}.\ea
\eqno(A.12)$$
This shows that the restriction of $\bar u(\cd\,;t,x)$ on $[\t,T]$
is not necessarily optimal for Problem (C) with initial pair
$(\t,\bar X(\t;t,x))$. Hence, Problem (C) is time-inconsistent.

\vfil\eject

Consider
\bel{}\left\{\ba{ll}
\ns\ds\dot X(s)=u(s),\qq s\in[t,T],\\
\ns\ds X(t)=x,\ea\right.\ee
with cost functional
\bel{}J(t,x;u(\cd))=\int_t^Tu(s)^2ds+h(t)X(T;t,x,u(\cd))^2,\ee
Let $\D:0=t_0<t_1<t_2<\cds<t_{N-1}<t_N=T$ be a partition of $[0,T]$.
Consider an LQ problem on $[t_{N-1},t_N]$, with the state equation
\bel{}\left\{\ba{ll}
\ns\ds\dot X_N(s)=u_N(s),\qq s\in[t_{N-1},t_N],\\
\ns\ds X_N(t_{N-1})=x,\ea\right.\ee
and cost functional
\bel{}J_N(t_{N-1},x;u_N(\cd))=\int_{t_{N-1}}^{t_N}u_N(s)^2ds+h(t_{N-1})X_N(t_N)^2.\ee
The corresponding Riccati differential equation reads
$$\left\{\ba{ll}
\ns\ds\dot P^N(s)-P^N(s)^2=0,\qq s\in[t_{N-1},t_N],\\
\ns\ds P^N(t_N)=h(t_{N-1}),\ea\right.
\eqno(A.1)$$
Simple calculation shows that
$$P^N(s)={P^N(t_N)\over1+P^N(t_N)(t_N-s)}={h(t_{N-1})\over1+h(t_{N-1})(t_N-s)},\qq s\in[t_{N-1},t_N].
\eqno(A.2)$$
The optimal control trajectory is the solution to the following
closed-loop system
$$\left\{\ba{ll}
\ns\ds\dot{\bar X}_N(s)=-P^N(s)\bar X_N(s),\qq s\in[t_{N-1},t_N],\\
\ns\ds\bar X_N(t_{N-1})=x,\ea\right.
\eqno(A.3)$$
which is given by
$$\bar X_N(s)=x{1+P^N(t_N)(t_N-s)\over1+P^N(t_N)(t_N-t_{N-1})}=x{1+h(t_{N-1})(t_N-s)\over1+h(t_{N-1})(t_N-t_{N-1})}\,,\qq
s\in[t_{N-1},t_N],
\eqno(A.4)$$
and the optimal control is given by
$$\bar u_N(s)=-P^N(s)\bar
X_N(s)=-{xP^N(t_N)\over1+P^N(t_N)(t_N-t_{N-1})}=-{xh(t_{N-1})\over1+h(t_{N-1})(t_N-t_{N-1})},\qq
s\in[t,T].
\eqno(A.5)$$
Now, on $[t_{N-2},t_{N-1}]$, we consider state equation
\bel{}\left\{\ba{ll}
\ns\ds\dot X_{N-1}(s)=u_{N-1}(s),\qq s\in[t_{N-2},t_{N-1}],\\
\ns\ds X_{N-1}(t_{N-2})=x,\ea\right.\ee
with the cost functional
$$\ba{ll}
\ns\ds
J_{N-1}(t_{N-2},x;u_{N-1}(\cd))=\int_{t_{N-2}}^{t_{N-1}}u_{N-1}(s)^2ds+\int_{t_{N-1}}^{t_N}\bar u_N(s)^2ds
+h(t_{N-2})\bar X(t_N)^2\\
\ns\ds=\int_{t_{N-2}}^{t_{N-1}}u_{N-1}(s)^2ds+{h(t_{N-1})^2(t_N-t_{N-1})+h(t_{N-2})\over[1+h(t_{N-1})(t_N-t_{N-1})]^2}
X(t_{N-1})^2\\
\ns\ds=\int_{t_{N-2}}^{t_{N-1}}u_{N-1}(s)^2ds+{P^N(t_N)^2(t_N-t_{N-1})+h(t_{N-2})\over[1+P^N(t_N)
(t_N-t_{N-1})]^2}
X(t_{N-1})^2\ea$$
For the corresponding LQ problem, the Riccati equation is
$$\left\{\ba{ll}
\ns\ds\dot P^{N-1}(s)-P^{N-1}(s)^2=0,\qq s\in[t_{N-2},t_{N-1}],\\
\ns\ds
P^{N-1}(t_{N-1})={P^N(t_N)^2(t_N-t_{N-1})+h(t_{N-2})\over[1+P^N(t_N)(t_N-t_{N-1})]^2}.\ea\right.
\eqno(A.1)$$
Simple calculation shows that
$$\ba{ll}
\ns\ds P^{N-1}(s)={P^{N-1}(t_{N-1})\over1+P^{N-1}(t_{N-1})(t_{N-1}-s)}\\
\ns\ds={{P^N(t_N)^2(t_N-t_{N-1})+h(t_{N-2})\over[1+P^N(t_N)(t_N-t_{N-1})]^2}
\over1+{P^N(t_N)^2(t_N-t_{N-1})+h(t_{N-2})\over[1+P^N(t_N)(t_N-t_{N-1})]^2}(t_{N-1}-s)}\\
\ns\ds={P^N(t_N)^2(t_N-t_{N-1})+h(t_{N-2})
\over[1+P^N(t_N)(t_N-t_{N-1})]^2+[P^N(t_N)^2(t_N-t_{N-1})+h(t_{N-2})](t_{N-1}-s)},\\
\ns\ds={h(t_{N-1})^2(t_N-t_{N-1})+h(t_{N-2})
\over[1+h(t_{N-1})(t_N-t_{N-1})]^2+[h(t_{N-1})^2(t_N-t_{N-1})+h(t_{N-2})](t_{N-1}-s)},\\
\ns\ds\qq\qq\qq\qq\qq\qq\qq\qq s\in[t_{N-2},t_{N-1}].\ea
\eqno(A.2)$$
Note
$$\ba{ll}
\ns\ds
P^N(t_{N-1})-P^{N-1}(t_{N-1})\\
\ns\ds={P^N(t_N)\over1+P^N(t_N)(t_N-t_{N-1})}-{P^N(t_N)^2(t_N-t_{N-1})+h(t_{N-2})\over[1+P^N(t_N)
(t_N-t_{N-1})]^2}\\
\ns\ds={P^N(t_N)[1+P^N(t_N)(t_N-t_{N-1})]-P^N(t_N)^2(t_N-t_{N-1})-h(t_{N-2})\over[1+P^N(t_N)
(t_N-t_{N-1})]^2}\\
\ns\ds={P^N(t_N)-h(t_{N-2})\over[1+P^N(t_N)(t_N-t_{N-1})]^2}
={h(t_{N-1})-h(t_{N-2})\over[1+h(t_{N-1})(t_N-t_{N-1})]^2}\ea$$
The optimal trajectory is the solution to the following
closed-loop system
$$\left\{\ba{ll}
\ns\ds\dot{\bar X}_{N-1}(s)=-P^{N-1}(s)\bar X_{N-1}(s),\qq s\in[t_{N-2},t_{N-1}],\\
\ns\ds\bar X_{N-1}(t_{N-2})=x,\ea\right.
\eqno(A.3)$$
which is given by
$$\bar X_{N-1}(s)=x{1+P^{N-1}(t_{N-1})(t_{N-1}-s)\over1+P^{N-1}(t_{N-1})(t_{N-1}-t_{N-2})}\,,\qq
s\in[t_{N-2},t_{N-1}],
\eqno(A.4)$$
and the optimal control is given by
$$\bar u_{N-1}(s)=-P^{N-1}(s)\bar
X_{N-1}(s)=-{xP^{N-1}(t_{N-1})\over1+P^{N-1}(t_{N-1})(t_{N-1}-t_{N-2})},\qq
s\in[t,T].
\eqno(A.5)$$

\ms

Now, on $[t_{N-3},t_{N-2}]$, we consider state equation
\bel{}\left\{\ba{ll}
\ns\ds\dot X_{N-2}(s)=u_{N-2}(s),\qq s\in[t_{N-3},t_{N-2}],\\
\ns\ds X_{N-2}(t_{N-3})=x,\ea\right.\ee
with the cost functional
$$\ba{ll}
\ns\ds
J_{N-2}(t_{N-3},x;u_{N-2}(\cd))\\
\ns\ds=\int_{t_{N-3}}^{t_{N-2}}u_{N-2}(s)^2ds+\int_{t_{N-2}}^{t_{N-1}}\bar u_{N-1}(s)^2ds
+\int_{t_{N-1}}^{t_N}\bar u_N(s)^2ds+h(t_{N-3})\bar X(t_N)^2\\
\ns\ds=\int_{t_{N-3}}^{t_{N-2}}u_{N-2}(s)^2ds+{P^{N-1}(t_{N-1})^2
(t_{N-1}-t_{N-2})\over[1+P^{N-1}(t_{N-1})(t_{N-1}-t_{N-2})]^2}X_{N-2}(t_{N-2})^2\\
\ns\ds\qq\qq+{P^N(t_N)^2(t_N-t_{N-1})+h(t_{N-3})\over[1+P^N(t_N)(t_N-t_{N-1})]^2}\bar X_{N-1}(t_{N-1})^2
\\
\ns\ds=\int_{t_{N-3}}^{t_{N-2}}u_{N-2}(s)^2ds+{P^{N-1}(t_{N-1})^2
(t_{N-1}-t_{N-2})\over[1+P^{N-1}(t_{N-1})(t_{N-1}-t_{N-2})]^2}X_{N-2}(t_{N-2})^2\\
\ns\ds\qq\qq+{P^N(t_N)^2(t_N-t_{N-1})+h(t_{N-3})\over[1+P^N(t_N)(t_N-t_{N-1})]^2
[1+P^{N-1}(t_{N-1})(t_{N-1}-t_{N-2})]^2}X_{N-2}(t_{N-2})^2
\\
\ns\ds=\int_{t_{N-2}}^{t_{N-1}}u_{N-1}(s)^2ds+{P^N(t_N)^2(t_N-t_{N-1})+h(t_{N-2})\over[1+P^N(t_N)
(t_N-t_{N-1})]^2}
X(t_{N-1})^2\ea$$
For the corresponding LQ problem, the Riccati equation is
$$\left\{\ba{ll}
\ns\ds\dot P^{N-1}(s)-P^{N-1}(s)^2=0,\qq s\in[t_{N-2},t_{N-1}],\\
\ns\ds
P^{N-1}(t_{N-1})={P^N(t_N)^2(t_N-t_{N-1})+h(t_{N-2})\over[1+P^N(t_N)(t_N-t_{N-1})]^2}.\ea\right.
\eqno(A.1)$$
Simple calculation shows that
$$\ba{ll}
\ns\ds P^{N-1}(s)={P^{N-1}(t_{N-1})\over1+P^{N-1}(t_{N-1})(t_{N-1}-s)}\\
\ns\ds={{P^N(t_N)^2(t_N-t_{N-1})+h(t_{N-2})\over[1+P^N(t_N)(t_N-t_{N-1})]^2}
\over1+{P^N(t_N)^2(t_N-t_{N-1})+h(t_{N-2})\over[1+P^N(t_N)(t_N-t_{N-1})]^2}(t_{N-1}-s)}\\
\ns\ds={P^N(t_N)^2(t_N-t_{N-1})+h(t_{N-2})
\over[1+P^N(t_N)(t_N-t_{N-1})]^2+[P^N(t_N)^2(t_N-t_{N-1})+h(t_{N-2})](t_{N-1}-s)},\\
\ns\ds={h(t_{N-1})^2(t_N-t_{N-1})+h(t_{N-2})
\over[1+h(t_{N-1})(t_N-t_{N-1})]^2+[h(t_{N-1})^2(t_N-t_{N-1})+h(t_{N-2})](t_{N-1}-s)},\\
\ns\ds\qq\qq\qq\qq\qq\qq\qq\qq s\in[t_{N-2},t_{N-1}].\ea
\eqno(A.2)$$
Note
$$\ba{ll}
\ns\ds
P^N(t_{N-1})-P^{N-1}(t_{N-1})\\
\ns\ds={P^N(t_N)\over1+P^N(t_N)(t_N-t_{N-1})}-{P^N(t_N)^2(t_N-t_{N-1})+h(t_{N-2})\over[1+P^N(t_N)
(t_N-t_{N-1})]^2}\\
\ns\ds={P^N(t_N)[1+P^N(t_N)(t_N-t_{N-1})]-P^N(t_N)^2(t_N-t_{N-1})-h(t_{N-2})\over[1+P^N(t_N)
(t_N-t_{N-1})]^2}\\
\ns\ds={P^N(t_N)-h(t_{N-2})\over[1+P^N(t_N)(t_N-t_{N-1})]^2}
={h(t_{N-1})-h(t_{N-2})\over[1+h(t_{N-1})(t_N-t_{N-1})]^2}\ea$$
The optimal trajectory is the solution to the following
closed-loop system
$$\left\{\ba{ll}
\ns\ds\dot{\bar X}_{N-1}(s)=-P^{N-1}(s)\bar X_{N-1}(s),\qq s\in[t_{N-2},t_{N-1}],\\
\ns\ds\bar X_{N-1}(t_{N-2})=x,\ea\right.
\eqno(A.3)$$
which is given by
$$\bar X_{N-1}(s)=x{1+P^{N-1}(t_{N-1})(t_{N-1}-s)\over1+P^{N-1}(t_{N-1})(t_{N-1}-t_{N-2})}\,,\qq
s\in[t_{N-2},t_{N-1}],
\eqno(A.4)$$
and the optimal control is given by
$$\bar u_{N-1}(s)=-P^{N-1}(s)\bar
X_{N-1}(s)=-{xP^{N-1}(t_{N-1})\over1+P^{N-1}(t_{N-1})(t_N-t_{N-1})},\qq
s\in[t,T].
\eqno(A.5)$$

\vfil\eject

 Now, if we let
$$J(t;\t,y;u(\cd))=\int_\t^Tu(s)^2ds+h(t)X(T;\t,y,u(\cd))^2,\qq
\t\in[t,T],
\eqno(A.6)$$
then the optimal value function (for fixed $t$) is given by
$$\ba{ll}
\ns\ds V(t;\t,y)\equiv\inf_{u(\cd)\in\cU[\t,T]}
J(t;\t,y;u(\cd))=J(t;\t,y;\bar u(\cd))=P(\t;t)y^2\\
\ns\ds\qq\q~~={h(t)\over
1+h(t)(T-\t)}y^2,\qq\forall(\t,y)\in[t,T]\times\dbR.\ea
\eqno(A.7)$$
Next, let $\t\in(t,T)$, we consider Problem (C) on $[\t,T]$ with
initial state
$$y=\bar X(\t;t,x)=x{1+h(t)(T-\t)\over1+h(t)(T-t)}.
\eqno(A.8)$$
The same as above, we see that the corresponding solution to the
Riccati equation is given by
$$P(s;\t)={h(\t)\over1+h(\t)(T-s)},\qq s\in[\t,T],
\eqno(A.9)$$
and
$$\ba{ll}
\ns\ds\inf_{u(\cd)\in\cU[\t,T]}J(\t;\t,y;u(\cd))=P(\t;\t)y^2={h(\t)y^2\over
1+h(\t)(T-\t)}\\
\ns\ds={x^2h(\t)[1+h(t)(T-\t)]^2\over[1+h(\t)(T-\t)][1+h(t)(T-t)]^2}.\ea
\eqno(A.10)$$
However,
$$\ba{ll}
\ns\ds J(\t,y;\bar u(\cd))=\int_\t^T\bar u(s)^2ds+h(\t)
X(T;\t,y,\bar u(\cd))^2\\
\ns\ds={x^2h(t)^2(T-\t)\over[1+h(t)(T-t)]^2}+h(\t)\[y-{xh(t)(T-\t)\over1+h(t)(T-t)}\]^2\\
\ns\ds={x^2h(t)^2(T-\t)\over[1+h(t)(T-t)]^2}+{x^2h(\t)\over[1+h(t)(T-t)]^2}.\ea
\eqno(A.11)$$
Hence,
$$\ba{ll}
\ns\ds J(\t,y;\bar
u(\cd))-\inf_{u(\cd)\in\cU[\t,T]}J(\t,y;u(\cd))\\
\ns\ds={x^2h(t)^2(T-\t)+x^2h(\t)\over[1+h(t)(T-t)]^2}-{x^2h(\t)[1+h(t)(T-\t)]^2\over[1+h(\t)(T-\t)]
[1+h(t)(T-t)]^2}\\
\ns\ds=x^2\[{[h(t)^2(T-\t)+h(\t)][1+h(\t)(T-\t)]-h(\t)[1+h(t)(T-\t)]^2\over[1+h(\t)(T-\t)][1+h(\t)(T-t)]^2}\]\\
\ns\ds={x^2[h(\t)-h(t)]^2(T-\t)\over[1+h(\t)(T-t)][1+h(t)(T-\t)]^2}>0,\qq\qq\hb{unless
$h(\t)=h(t)$ or $x=0$}.\ea
\eqno(A.12)$$
This shows that the restriction of $\bar u(\cd\,;t,x)$ on $[\t,T]$
is not necessarily optimal for Problem (C) with initial pair
$(\t,\bar X(\t;t,x))$. Hence, Problem (C) is time-inconsistent.

\end{document}